\documentclass[12pt]{article}

\usepackage{amssymb}
\usepackage{amsmath}

\def \msk {\medskip}
\def \bsk {\bigskip}

\begin{document}

\begin{center}
{\large Precoloring Extension.\ \,IV.\, General Bounds and List Colorings}

\bsk

Mih\'{a}ly Hujter; Zsolt Tuza
\end{center}

\bsk\bsk

\sf

{\footnotesize
\begin{quote}
 \textbf{Abstract.} A guaranteed upper bound is proved for the
  time complexity of the list-coloring problem on graphs.
\end{quote}
}

\bsk

For a positive integer $n$, the set $\{1,\ldots ,n\}$ is denoted by $%
[1,n]$. 
Throughout this note, $n$ will denote the \emph{order} of a graph,
that is the number of vertices.

A \emph{proper} $k$-\emph{coloring} (or simply, a $k$-coloring) of a
graph $G=(V,E)$ is a function $f:V\rightarrow \lbrack 1,k]$ such that $%
f(v)\neq f(w)$ whenever $vw\in E$. The \emph{vertex coloring problem}, whose
input is a graph $G$ and a positive integer $k$, consists of deciding
whether $G$ is properly $k$-colorable or not. In the \emph{list-coloring
problem}, every vertex $v$ comes equipped with a list of permitted colors $%
L(v)\subseteq \lbrack 1,\kappa ]$ for some given positive integer $\kappa $,
and we require the proper coloring to respect these lists, i.e., $f(v)\in
L(v)$. 

We say that a constant real number $t>1$ \emph{fits} $G$ if for any $m$%
-element subset $M$ of $V$, the \emph{induced subgraph} $G[M]$ has at most $%
t^{m}$ \emph{inlcusion-maximal independent sets}. It is well-known that $%
1.44225$ fits any graph because $1.44225^{3}>3$; see \cite{M} for a proof of
the fact that an $n$-vertex graph contains at most $3^{n/3}$
inclusion-maximal independent sets. Similarly, $1.41422$ fits any
triangle-free graph because $1.41422^{2}>2$; see \cite{SIAM} for a proof of
the fact than a triangle-free graph contains at most $2^{n/2}$
inclusion-maximal independent sets.
Furthermore it was proved in \cite{Netw} that, for any fixed natural number $k$,
graphs not containing an induced matching with $k+1$ edges have at most
$O(n^{2k})$ inclusion-maximal independent sets.
Hence for every $k$ and every real $\epsilon>0$ there exists an 
$n_0=n_0(k,\epsilon)$ such that $t=1+\epsilon$ fits every such graph
of order $n>n_0$.

A special case of the list-coloring problem is the \emph{precoloring
extension} with $\kappa $ colors. Here each list is either $[1,\kappa ]$ or a
one-element subset of it. (See \cite{prext1,prext2,prext3} for more details
on precoloring extension.) Another special case is the $k$-\emph{choosability} 
problem where each color list has size $k$ for a fixed $k$. (See \cite{ERT,t97}
for more details on choosability.)

The main result in the present note, Theorem 1, claims that there exists a 
polynomial $p(n)$ such that if $t$ fits a graph $G$ with $n$ vertices,
then it can be decided in $p(n)\cdot (1+t)^{n}+O(\kappa n)$ time whether $G$ is
properly list-colorable or not. (Here the polynomial $p$ is independent of $%
t $, and the $O(\kappa n)$ term is irrelevant except for reading an input with
rather long vertex lists.) In case of a positive answer, we will construct a proper
list-coloring needing no extra time. We will extend the methods and results
of Lawler \cite{L}.

Given a graph $G=(V,E)$ on $n$ vertices, we fix a
permutation of the nonempty subsets of $V$ for which $W,U\subseteq V$ and $%
|W|>|U|$\ imply that $W$ is before $U$ in the fixed permutation. For any
nonempty $W\subseteq V$ let $\langle W\rangle $ denote the position of $W$
in the fixed permutation. As a special case, $\langle V\rangle =1$, and $%
|W|=n-1$ implies that $\langle W\rangle \in \{2,3,...,n+1\}$, furthermore $%
|W|=1\ $implies that $\langle W\rangle \in \{2^{n}-n,2^{n}-n+1,...,2^{n}-1\}$%
. We will consider a zero-one sequence $(a_{m}),\ m=1,2,...,2^{n}-1$. We say
that the sequence has the \emph{partial increasing property} (PIP, for
schort) if $W\supset U$ implies $a_{\langle W\rangle }\leq a_{\langle
U\rangle }$. In words, if $W$ is before $U$ in the permutation and if $%
a_{\langle W\rangle }=1$, then $a_{\langle U\rangle }=1$ holds by pip.

Initially, the sequence $(a_{m})$, $m=1,2,...,2^{n}-1$ will be set as the 
\emph{characteristic sequence of the independent sets} in $G$. In other
words, $a_{\langle W\rangle }=1$ holds if and only if the induced subgraph $%
G[W]$ is edgeless. Later some zero values in the sequence will be changed to
one; however, the pip will be managed.

Clearly, list-colorability can be reformulated as follows: The vertex
set $V$ can be partitioned into independent set \emph{color classes} $%
V_{1},\ldots ,V_{\kappa }$ such that for each $j\in \lbrack 1,\kappa ]$ and
for each $v\in V_{j}$ the relation $j\in L(v)$ holds. Here we allow the
empty set to occur among $V_{1},\ldots ,V_{\kappa }$ any number of times.

For each $j\in \lbrack 1,\kappa ]$ and for each $v\in V$ we define $L^{j}(v)$
as $L(v)\cap \lbrack 1,j]$. Clearly $L^{1}(v)\subseteq \ldots \subseteq
L^{\kappa }(v)=L(v)$.

For each $j\in \lbrack 1,\kappa ]$ we will consider the $j^{\text{th}}$
version of the list-coloring problem, denoted by $\mathcal{LC}^{j}$ for all
induced subgraphs $G[W]$ by considering the lists $L^{j}(v)$. Note that $%
\mathcal{LC}^{1}$ is trivial since $G[W]$ is list-colorable (with respect to
the lists $L^{j}(v)$, $w\in W$) if and only if, on the one hand, $W$ is an
independent set in $G$, and on the other hand, $1\in $ $L(w)$ for all $w\in
W $. Starting from the list $\mathcal{I}$ of all inclusion-maximal
independent sets of $G=(V,E)$, the sequence $(a_{m}),$ $m=1,2,...,2^{n}-1$
can be computed as follows:
\begin{enumerate}
 \item Let all $a_{m}$ get the value $0$.
 \item For each $I\in \mathcal{I}$, consider 
\[
I^{(1)}=\{v\in I:1\in L(v)\} 
\]%
and set $a_{\langle I^{(1)}\rangle }=1$.
 \item Manage the above mentioned PIP of the sequence $(a_{m})$.
\end{enumerate}
For this purpose we go forward in
the sequence and for any $a_{\langle W\rangle }=1$ and for any $w\in W$ we
set $a_{\langle W-w\rangle }=1.$ (Here $W-w$ simply denotes $W\diagdown
\{w\} $.)

\msk

\textbf{Lemma 1}. After performing the above three steps, the final values
in the sequence $(a_{\langle W\rangle })$ give the true answer for the
list-coloring problem for each $G[W]$ with respect to the lists $L^{1}(v)$.
This means that $a_{\langle W\rangle }=1$ holds if and only if $G[W]$ is
list-colorable with respect to the lists $L^{1}(v)$, $v\in W$.

\emph{Proof. }For any inclusion-maximal vertex set $I$, in step 2 the
sequence element $a_{\langle I\rangle }$ gets value $1$, and in step 3, for
any other independent set $W$, the final value of $a_{\langle W\rangle }$
will also be $1$. $\blacksquare $

\msk

Now for any fixed $j=2,3,...,\kappa $, we will recursively solve the
list-colorabilty problem $\mathcal{LC}^{j}\ $for all $G[W]$ by using the
lists consisting of the sets%
\[
I^{(j)}=\{v\in I:j\in L(v)\} 
\]%
obtained from all inclusion-maximal independent sets $I$ in $G$. For a fixed $%
j\in \{2,3,...,\kappa \}$, we assume that the problem $\mathcal{LC}^{j-1}$
is completely solved, and the result is properly recorded in the sequence $%
(a_{m}),\ m=1,2,...,2^{n}-1$, that is for any $G[W]$ this induced subgraph
is list-colorable with respect to the color lists $L^{j-1}(v)$, $v\in W$, if
and only if $a_{\langle W\rangle }=1$.

\msk

\textbf{Lemma 2.} The graph $G=(V,E)$ is list-colorable with respect to the color
lists $L^{j}(v)$, $v\in V$, if and only if there is at least one independent
set $I^{(j)}$ for which either $I^{(j)}=V$ or $a_{\langle V\diagdown
I^{(j)}\rangle }=1$ holds as a result of the $\mathcal{LC}^{j-1}$ problem.

\emph{Proof.} If $I^{(j)}=V$, then $G$ is clearly list-colorable since $%
f(v)=j$ suffices for all $v\in V$. If there is an $I^{(j)}\neq V$ for which $%
a_{\langle V\diagdown I^{(j)}\rangle }=1$ holds, then any list-coloring $f$
of $G[V\diagdown I^{(j)}]$ with respect to the lists $L^{j-1}(v)$, $v\in
V\diagdown I^{(j)}$, is clearly extendible to obtain a list-coloring of $G$
by defining $f(v)=j$ for all $v\in I^{(j)}$. On the other hand, in any
list-coloring $f$ of $G$ with respect to the lists $L^{j}(v),\ v\in V$, the
vertices $v$ with $f(v)=j$, must form an independent set $J$, and this $J$
must be contained in at least one $I^{(j)}$, and now the restriction of $f$
onto $V\diagdown I^{(j)}$ gives us a list-coloring of $G[V\diagdown I^{(j)}]$
with respect to the lists $L^{j-1}(v)$, $v\in V\diagdown I^{(j)}$. $%
\blacksquare $

\msk

Now Lemma 2 shows how we can solve the $\mathcal{LC}^{j}$ problem for $G$.
Similarly we can solve the $\mathcal{LC}^{j}$ problem for any induced
subgraph $G[W]$. Here we can use the same sequence $(a_{m}),\
m=1,2,...,2^{n}-1$. 
For each $W$ we scan the sets $I^{(j)}$ for all inclusion-maximal independent
sets $I$ in $G[W]$, and update $a_{\langle W\rangle }$ from $0$ to $1$ if
$a_{\langle W\diagdown I^{(j)}\rangle }=1$ holds in the sequence after round $j-1$.
In this procedure PIP ensures that the elements $a_{\langle U \rangle }$ may be 
updated only when they are not needed anymore for updates of any other elements.

If $t$ fits $G$, then the total number of steps is not
larger than $p(n)$ times the following number, where $p(n)$ is a polynomial
independent of $G$ and independent of $t$, too: 
\[
\sum_{W\neq \emptyset }\left( _{|W|}^{\ n}\right) t^{|W|} 
\]%
By the binomial theorem,$\ $the above sum is $(1+t)^{n}-1$. In
summary, we obtain the following theorem.

\msk

\textbf{Theorem 1.} There exists a polynomial $p$ for which if $t$ fits 
$G$, an $n$-vertex graph, then all the problems $\mathcal{LC}^{j}$, $%
j=1,2,...,\kappa $ can be solved in 
\[
 p(n)\cdot (1+t)^{n} + O(\kappa n)
\]%
time.

\emph{Proof.}
While reading the input, we put aside all vertices whose lists have length
at least $n$.
Once the rest of the graph is properly list-colored, the vertices with long lists
can be colored in an obvious greedy way.
For this reason we can assume from now on that $\kappa < n$ holds.

It is proved in \cite{TIAS} that if a graph on $n$ vertices has $N$
inclusion-maximal independent sets, then those $N$ sets can be listed in
$O(n^3 N)$ time.
We apply this method $\kappa$ times for all induced subgraphs $G[W]$.
The input graph $G$ is list colorable if and only if $a_{\langle V \rangle} = 1$
holds after round $\kappa$.
Thus there exists a suitable polynomial $p(n)=O(n^4)$.
$\blacksquare $

\msk

If, in case of list-colorability, we want to output a proper list-coloring $%
f:V\rightarrow \lbrack 1,\kappa ]$, we simply have to memorize an
independent set $I^{(\kappa )}$ for which $G[V\diagdown I^{(k)}]$ is
list-colorable with respect to the lists $L^{\kappa -1}(v)$, $v\in
V\diagdown I^{(\kappa )}$, and for any $w\in I^{(\kappa )}$ we can set $%
f(w)=\kappa $. Then we decrease the value of $\kappa $ by one, and we repeat
the procedure for the decreased graph instead of the old one, and for the
new $\mathcal{LC}^{\kappa }$ problem instead of the old $\mathcal{LC}%
^{\kappa }$ problem. And so on.


\bigskip

\begin{thebibliography}{9}
\bibitem{L} Lawler, E.\ L.: \emph{A note on the complexity of the chromatic
number problem}. Inf.\ Process.\ Lett.\  \textbf{5} (1976) 66--67.

\bibitem{TIAS} Tsukiyama, S.; Ide, M.; Ariyoshi, H.; Shirakawa, I.: \emph{A
new algorithm for generating all the maximal independent sets}. SIAM J.\
Comput.\  \textbf{6} (1977) 505--517.

\bibitem{M} Moon, J.; Moser, L.: \emph{On cliques in graphs}, Israel J.
Math. \textbf{3} (1965) 23--28.

\bibitem{ERT} Erd\H{o}s, P.; Rubin, A.L.; Taylor, H.: \emph{Choosability in
graphs}. Combinatorics, graph theory and computing, Proc. West Coast Conf.,
Arcata/Calif. 1979, (1980) 125--157.

\bibitem{prext1} Bir\'{o}, M.; Hujter, M.; Tuza, Zs.: \emph{Precoloring
extension.\ I.\ Interval graphs}. Discrete Math.\  \textbf{100} (1992)
267--279.

\bibitem{prext2} Hujter, M.; Tuza, Zs.: \emph{Precoloring extension.\ II.\
Graph classes related to bipartite graphs}. Acta Math.\ Univ.\ Comenianae 
\textbf{62} (1993) 1--11.

\bibitem{prext3} Hujter, M., and Tuza, Zs.: \emph{Precoloring extension.\
III.\ Classes of perfect graphs}. Combin.\ Probab.\ Comput.\  \textbf{5} (1996)
35--56.

\bibitem{SIAM} Hujter, M.; Tuza, Zs.: \emph{The number of maximal
independent sets in triangle-free graphs}. SIAM J.\ on Discrete
Math.\ \textbf{6} (1993) 284--288.

\bibitem{Netw} Farber, M.; Hujter, M.; Tuza, Zs.: \emph{An upper bound on
the number of cliques in a graph}. Networks \textbf{23} (1993) 207--210.

\bibitem{t97} Tuza, Zs.: \emph{Graph colorings with local constraints -- A survey}.
Discuss.\ Math.\ Graph Theory \textbf{17} (1993) 161--228.

\end{thebibliography}
\end{document}